\documentclass[12pt]{article}
\usepackage{graphicx}
\usepackage{color}
\catcode`\@=11 \@addtoreset{equation}{section}

\catcode`\@=12
\usepackage{amssymb}
\usepackage{amsmath}
\usepackage{smaller}

\usepackage{color}
\definecolor{light}{gray}{0.50}
\definecolor{heavy}{gray}{0.35}
\definecolor{black}{gray}{0.0}
\definecolor{dgreen}{rgb}{0.0,0.7,0}
\definecolor{dred}{rgb}{0.9959,0,0}
\definecolor{green}{rgb}{0.0,0.99599,0.0}
\definecolor{purple}{rgb}{0.6,0.0,0.4}

\newcounter{ewacommentno}
\newcounter{eduardcommentno}
\newcounter{antonincommentno}
\newcommand{\rupert}[1]{\textcolor{dgreen}{#1}}
\newcounter{rupertcommentno}
\newcommand{\eq}[1]{(\ref{#1})}
\newcommand{\dss}{\displaystyle}

\newtheorem{Theorem}{Theorem}[section]
\newtheorem{Proposition}{Proposition}[section]
\newtheorem{Lemma}{Lemma}[section]
\newtheorem{Corollary}{Corollary}[section]
\newtheorem{Remark}{Remark}[section]
\newtheorem{Definition}{Definition}[section]

\newcommand{\bTheorem}[1]{
\begin{Theorem} \label{T#1} }
\newcommand{\eT}{\end{Theorem}}

\newcommand{\bProposition}[1]{
\begin{Proposition} \label{P#1}}
\newcommand{\eP}{\end{Proposition}}

\newcommand{\bLemma}[1]{
\begin{Lemma} \label{L#1} }
\newcommand{\eL}{\end{Lemma}}

\newcommand{\bCorollary}[1]{
\begin{Corollary} \label{C#1} }
\newcommand{\eC}{\end{Corollary}}

\newcommand{\bRemark}[1]{
\begin{Remark} \label{R#1} }
\newcommand{\eR}{\end{Remark}}

\newcommand{\bDefinition}[1]{
\begin{Definition} \label{D#1} }
\newcommand{\eD}{\end{Definition}}

\newcommand{\bFormula}[1]{
\begin{equation} \label{#1}}
\newcommand{\eF}{\end{equation}}

\newcommand{\The}{\Theta_\ep}

\newcommand{\Ov}[1]{\overline{#1}}

\newcommand{\DC}{C^\infty_c}

\newcommand{\vr}{\varrho}
\newcommand{\vre}{\vr_\ep}

\newcommand{\vue}{\vu_\ep}

\newcommand{\vu}{\vc{u}}
\newcommand{\vc}[1]{{\boldsymbol{#1}}}
\newcommand{\qed}{\bigskip \rightline {Q.E.D.} \bigskip}
\newcommand{\Div}{{\rm div}_x}
\newcommand{\Grad}{\nabla_x}

\newcommand{\tn}[1]{\mathbb{#1}}
\newcommand{\dx}{{\rm d} {x}}
\newcommand{\dt}{{\rm d} t }

\newcommand{\dxdt}{\dx \ \dt}
\newcommand{\intO}[1]{\int_{\Omega} #1 \ \dx}

\newcommand{\bProof}{{\bf Proof: }}

\newcommand{\mT}{\mathcal{T}}

\newcommand{\ep}{\varepsilon}

\newcommand{\R}{\mathbb{R}}

\definecolor{Cgrey}{rgb}{0.85,0.85,0.85}
\definecolor{Cblue}{rgb}{0.50,0.85,0.85}
\definecolor{Cred}{rgb}{1,0,0}
\definecolor{fancy}{rgb}{0.10,0.85,0.10}

\newcommand\Cbox[2]{%
    \newbox\contentbox%
    \newbox\bkgdbox%
    \setbox\contentbox\hbox to \hsize{%
        \vtop{
            \kern\columnsep
            \hbox to \hsize{%
                \kern\columnsep%
                \advance\hsize by -2\columnsep%
                \setlength{\textwidth}{\hsize}%
                \vbox{
                    \parskip=\baselineskip
                    \parindent=0bp
                    #2
                }%
                \kern\columnsep%
            }%
            \kern\columnsep%
        }%
    }%
    \setbox\bkgdbox\vbox{
        \color{#1}
        \hrule width  \wd\contentbox %
               height \ht\contentbox %
               depth  \dp\contentbox
        \color{black}
    }%
    \wd\bkgdbox=0bp%
    \vbox{\hbox to \hsize{\box\bkgdbox\box\contentbox}}%
    \vskip\baselineskip%
}

\date{}

\begin{document}


\title{On singular limits arising in the scale analysis of stratified fluid flows}

\author{Eduard Feireisl$^1$\thanks{The research of E.F.~leading to these results has received funding from the European Research Council under the European Union's Seventh Framework Programme (FP7/2007-2013)/ ERC Grant Agreement 320078. The Institute of Mathematics of the Academy of Sciences of the Czech Republic is supported by RVO:67985840.}  \and Rupert Klein$^2$\thanks{R.K.'s research has been partially funded by Deutsche Forschungsgemeinschaft (DFG) through grant CRC 1114.}\and Anton\' \i n Novotn\' y$^3$ \and Ewelina Zatorska$^{4,5}$\thanks{E.Z.'s research has been supported by the National Science Centre, Poland, grant 2014/14/M/ST1/00108 and by the fellowship START of the Foundation for Polish Science.}   }

\maketitle

\bigskip

\centerline{$^1$Institute of Mathematics of the Academy of Sciences of the Czech Republic}

\centerline{\v Zitn\' a 25, CZ-115 67 Praha 1, Czech Republic}

\medskip

\centerline{$^2$FB Mathematik \& Informatik, Freie Universit\"at Berlin}

\centerline{Arnimallee 6, 14195 Berlin}

\medskip

\centerline{$^3$IMATH, EA 2134, Université du Sud Toulon-Var}

\centerline{BP 20132, 83957 La Garde, France}

\medskip

\centerline{$^4$ Institute of Mathematics, Polish Academy of Sciences}

\centerline{ul. \'Sniadeckich 8, 00-656 Warszawa, Poland}

\medskip

\centerline{$^5$ Institute of Applied Mathematics and Mechanics, University of Warsaw}

\centerline{ul. Banacha 2, 02-097 Warszawa, Poland}






\maketitle

\bigskip





\begin{abstract}
We study the low Mach low Freude numbers limit in the compressible Navier-Stokes equations and the transport equation for evolution of an entropy variable -- the potential temperature $\Theta$. We consider the case of well-prepared initial data on "flat" tours and Reynolds number tending to infinity, and the case of ill-prepared data on an infinite slab. In both cases, we show that the weak solutions to the primitive system converge to the solution to the anelastic Navier-Stokes system and the transport equation for the second order variation of $\Theta$.

\end{abstract}


\tableofcontents

\section{Introduction}
\label{i}

{We study certain singular limits arising in the analysis of stratified fluid flows in meteorology and astrophysics, \cite{AlmgrenEtAl2006,Klein2010}. Moving beyond the case of homentropic flow, which was studied rigorously by Masmoudi, \cite{MAS2}, we include variations of an entropy variable -- the \emph{potential temperature}, $\Theta$, at second order in the Mach number. This is motivated by the observation that only when such entropy variations are accounted for, will the fluid flow equations support internal gravity waves, and these are responsible for a host of important physical processes, see, e.g., \cite{AlexanderEtAl2010}. Specifically, we study the asymptotic regime introduced originally in a formal asymptotic analysis by Ogura and Phillips \cite{OguraPhillips1962}. Within this regime, advection and internal waves act on comparable time scales while sound waves are asymptotically fast in the low Mach number limit. Thus, considering this flow regime allows us to incorporate the effects of entropy stratification without the added difficulty of an asymptotic three-scale problem, \cite{KleinEtAl2010}, that would have advection, internal waves, and sound act on three asymptotically separated time scales.} 

{An additional simplifying assumption concerns the time evolution of the potential temperature. In the mentioned application areas, some mean potential temperature stratification is generally observed over long time scales which differs from that associated with constant temperature, $T$. These stratifications are maintained against dissipative processes that tend to homogenize the temperature by various mechanisms including radiation, turbulent motions associated with convective instabilities, latent heat conversion, and similar ``diabatic effects''. In the present first approach to addressing flows with potential temperature stratification, we exclude the mathematical complications induced by these processes and instead ensure the maintainance of a mean potential temperature stratification by excluding diabatic processes and the influence of molecular transport on potential temperature, and by imposing a mean stratification in the initial data.}

{Under these premises, the fluid motion is described by the Navier-Stokes system in the isentropic regime. Specifically we analyze the following dimensionless system of equations for $\ep\ll 1$:}
\bFormula{i1}
\partial_t \vr + \Div (\vr \vu) = 0,
\eF
\bFormula{i2}
\partial_t (\vr \vu) + \Div (\vr \vu \otimes \vu) + \frac{1}{\ep^2} \Grad (\vr \Theta)^\gamma = \nu \Div \tn{S}(\Grad \vu) + \frac{1}{\ep^2} \vr \Grad F
\eF
\bFormula{i3}
\partial_t (\vr \Theta) + \Div (\vr \Theta \vu ) = 0,
\eF
where $\vr$ denotes the mass \emph{density}, $\vu$ the \emph{velocity}, and $\tn{S}(\Grad \vu)$ the viscous stress tensor, here given by Newton's rheological law \bFormula{i3a}
\tn{S}(\Grad \vu) = \mu \left( \Grad \vu + \Grad^t \vu - \frac{2}{3} \Div \vu\, \tn{I} \right) + \lambda \Div \vu\, \tn{I}, \ \mu > 0, \ \lambda \geq 0.
\eF
{Here we have adopted a distinguished limit that equates the Mach and (external wave) Froude numbers with the small parameter, i.e., ${\rm Ma} = {\rm Fr} = \ep$, while the Reynolds number ${\rm Re} = 1/\nu$ is treated as an independent parameter. In part of this work, we consider limit processes, in which $\ep\to 0$ and $\nu \to 0$ independently. The Strouhal number is unity as we adopt the advective time scale $t_{\rm ref} = \ell_{\rm ref}/u_{\rm ref}$ as the reference time scale for non-dimensionalization and, following the discussion of the last paragraph, we set the Prandtl number, ${\rm Pr} = \infty$, thereby neglecting the molecular transport of heat.}

The fluid occupies a slab
\bFormula{i3b}
\Omega = \Omega_h \times (0,1),
\eF
where the horizontal projection $\Omega_h$ is either a periodic ``flat'' torus $T_{2D} = \left( [0,1]|_{\{ 0,1 \}} \right)^2$ or the unbounded space
$\R^2$. The velocity is supposed to satisfy the complete slip {closed lid} boundary conditions,
\bFormula{i3c}
\vu \cdot \vc{n} = u_3|_{\partial \Omega} = 0, \ \left( \tn{S}(\Grad \vu) \cdot \vc{n} \right) \times \vc{n}|_{\partial \Omega} = 0,
\eF
where $\vc{n}$ denotes the outer normal vector to $\partial \Omega$. We shall denote $\vc{x} = [ \vc{x}_h,z]$, where $\vc{x}_h = [x_1,x_2]$

We take $F = - gz$, $g > 0$; whence $\Grad F = [0,0,-g]$ represents the effect of the \emph{gravitational force} acting on the fluid in the vertical direction. Accordingly, {and since we assume the potential temperature to be constant to leading order (see \eq{i4} below), the leading order} equilibrium distribution of the density $\tilde \vr$ satisfies
\bFormula{i3d}
\Grad {\tilde \vr}^\gamma = \tilde \vr \Grad F.
\eF
Obviously, equation (\ref{i3d}) admits a continuum of solutions, here we suppose $\tilde \vr = \tilde \vr (z)$ to be given, $\tilde \vr(z) > 0$ for all
$z \in [0,1]$.

Problem (\ref{i1}--\ref{i3c}) is supplemented with the initial data
\bFormula{i4}
\left\{ \begin{array}{c}
\vr(0,\cdot) = \vr_{0, \ep} = \tilde \vr + \ep \vr^{(1)}_{0, \ep},
\\ \\
\vu(0, \cdot) = \vu_{0,\ep},
\\ \\
\Theta(0, \ep) = \Theta_{0,\ep} = 1 + \ep^2 \Theta^{(2)}_{0,\ep}.
\end{array} \right\}
\eF
Moreover, if $\Omega_h = \R^2$ is unbounded, the far-field conditions
\bFormula{i4a}
\vr \to \tilde \vr, \ \Theta \to 1 \ \mbox{as}\ |x| \to \infty
\eF
are prescribed.


Our goal is to study the behavior of solutions $[\vre, \vue, \The]$ in the asymptotic limit $\ep \to 0$. Formally, it is not difficult to check that
\begin{equation}
\vr_\ep \approx \tilde \vr \ \mbox{as}\ \ep \to 0,
\end{equation}
{whereby} (\ref{i1}) reduces to the so-called \emph{anelastic constraint}
\bFormula{i5a}
\Div (\tilde \vr \vc{v}) = 0
\eF
for the limit velocity $\vc{v}$. Moreover, the specific form of the initial data for $\Theta$ indicates that
\begin{equation}
\Theta^{(2)}_\ep = \frac{ \The - 1 }{\ep^2} \approx \mT \ \mbox{for}\ \ep \to 0,
\end{equation}
where $\mT$ is transported by $\vc{v}$,
\bFormula{i5b}
\partial_t \mT + \vc{v} \cdot \Grad \mT = 0.
\eF
Finally, we also perform the limit in the momentum equation (\ref{i2}). Assuming that also the coefficient $\nu \to 0$ vanishes in the asymptotic
limit we recover, again formally, the system
\bFormula{i15}
\partial_t \vc{v} + \vc{v} \cdot \Grad \vc{v} + \Grad \Pi = - \mT \Grad F,
\eF
where $\Pi$ is the pressure, the presence of which is enforced by (\ref{i5a}). In the case $\nu > 0$ is kept fixed, a viscous tensor
$\frac{1}{\tilde \vr} \tn{S}(\Grad \vc{v})$ will necessarily appear in the limit system.





We consider {\emph{weak} solutions to the primitive system (\ref{i1}--\ref{i3}) according to the definition given in section~\ref{w}}. The advantage of such an approach is that these solutions are known to exist globally in time, the drawback is their low regularity given by the typically very poor {\it a priori} bounds. Two kinds of initial data will be considered: {\bf (i)} well-prepared data mimicking the {structure} of the data of the expected limit system, {\bf (ii)}  ill-prepared data {which merely require uniform bounds for certain rescaled quantities as $\ep \to 0$}. We show that in both cases the weak solutions $[\vre, \vue, \The]$ of the primitive system converge to the expected limit. In the case of the well-prepared data, we also obtain an explicit rate of decay in terms of $\ep$ and $\nu$ assumed to vanish in the asymptotic limit. 
The method used is that of relative energy (entropy) providing a suitable ``distance'' between the weak solutions of the primitive system and that of the target system. Accordingly, the convergence takes place on the life-span of (smooth) solutions to the limit problem.
Ill-prepared initial data give rise to high-frequency acoustic waves that must be filtered out in the limit. To this end, we consider the unbounded slab
$\Omega = \R^2 \times (0,1)$ allowing for dispersion of acoustic waves that actually vanish in the asymptotic limit on any \emph{compact} subset of $\Omega$. The analysis at this point relies on the application of the celebrated RAGE theorem. In contrast with the {case of} well-prepared data, we establish global-in-time convergence towards a weak solution of the limit problem. Unfortunately, as \emph{uniqueness} of the weak solutions for the limit problem is not known, there is no explicit
rate of decay available for the ill-prepared data.

The paper is organized as follows. In Section \ref{w}, we recall the concept of weak solution for both the primitive and the target system {used in this paper} and state the main results.
In Section \ref{R}, we introduce the relative (modulated) energy functional {needed in this context} and derive the necessary uniform bounds \emph{independent} of the scaling parameter $\ep$. Section \ref{A} is devoted to the case of well-prepared initial data while Section \ref{I} deals with ill-prepared data. Possible extensions and further comments on the methods used in the paper are discussed in Section \ref{C}.

\section{Weak solutions, main results}
\label{w}

In this section, we introduce the concept of weak solution to the Navier-Stokes system (\ref{i1}--\ref{i3}) and state our main result.

\subsection{Weak solutions to the primitive system}

We consider the weak solutions to system (\ref{i1}--\ref{i3}) belonging to the class
\bFormula{w1}
\begin{array}{rcl}
\vr 
  & \in 
    & C_{\rm weak}([0,T]; L^\gamma (\Omega)), 
      \\[10pt]
\Theta 
  & \in 
    & L^\infty((0,T) \times \Omega), 
    \\[10pt]
\vu 
  & \in 
    & L^2({0,T}; W^{1,2}(\Omega;\R^3)),
\end{array}
\eF
and enjoying further regularity and integrability properties allowed by the available {\it a priori} bounds. We remark that the impermeability constraint
\bFormula{w4}
\vu \cdot \vc{n}|_{\partial \Omega} = 0
\eF
makes sense in the above specified class.

\subsubsection{Weak formulation of the transport equations}

We say that $\vr$, $\vu$ is a weak solution to the equation of continuity (\ref{i1}) in $(0,T) \times \Omega$ if the integral identity
\bFormula{w4a}
\left[ \intO{ \vr \varphi } \right]_{t=\tau_1}^{t = \tau_2}
= \int_{\tau_1}^{\tau_2} \intO{ \left[ \vr \partial_t \varphi + \vr \vu \cdot \Grad \varphi \right] } \ \dt
\eF
holds for any $0 \leq \tau_1 \leq \tau_2\leq T$ and any test function $\varphi \in \DC([0,T] \times \Ov{\Omega})$.

As for equation (\ref{i3}), we consider its \emph{renormalized} version, specifically,
\bFormula{w5}
\left[ \intO{ \vr G(\Theta) \varphi } \right]_{t=\tau_1}^{t = \tau_2}
= \int_{\tau_1}^{\tau_2} \intO{ \left[ \vr G(\Theta) \partial_t \varphi + \vr G(\Theta) \vu \cdot \Grad \varphi \right] } \ \dt
\eF
for any $\varphi \in \DC([0,T] \times \Ov{\Omega})$ and any $G \in C(\rupert{\R})$.
\bRemark{w1}
Alternatively, we can replace (\ref{w5}) by postulating the renormalized version of the pure transport equation
\begin{equation}
\partial_t \Theta + \vu \cdot \Grad \Theta = 0,
\end{equation}
namely
\bFormula{w5a}
\left[ \intO{ G(\Theta) \varphi } \right]_{t=\tau_1}^{t = \tau_2}
= \int_{\tau_1}^{\tau_2} \intO{ \left[ G(\Theta) \partial_t \varphi + G(\Theta) \vu \cdot \Grad \varphi + G(\Theta) \Div \vu \varphi \right] } \ \dt
\eF
for any $\varphi \in \DC([0,T] \times \Ov{\Omega})$ and any $G \in C({\R})$.

\eR

\subsubsection{Weak formulation of the momentum balance, energy inequality}

The standard weak formulation of the momentum balance (\ref{i2}) reads:
\bFormula{w7}
\begin{array}{rcl}
\dss \left[ \intO{ \vr \vu \cdot \varphi } \right]_{t = \tau_1}^{t = \tau_2}
  & = 
    & \dss \int\limits_{\tau_1}^{\tau_2} 
      \int\limits_\Omega
      \left[ \vr \vu \cdot \partial_t \varphi + \vr \vu \otimes \vu : \Grad \varphi
      + \frac{1}{\ep^2} (\vr \Theta)^{\gamma} \Div \varphi\right. 
      \\[20pt] 
  & 
    & \dss \left.  
     \qquad\quad  - \ \nu \tn{S}(\Grad \vu) : \Grad \varphi 
      + \frac{1}{\ep^2} \vr \Grad F \cdot \varphi \right] \ \dt
\end{array}
\eF
for any $\varphi \in \DC([0,T] \times \Ov{\Omega}; \R^3)$, $\varphi \cdot \vc{n}|_{\partial \Omega} = 0$.

In addition, we focus on the class of \emph{finite energy} weak solutions satisfying the energy inequality
\bFormula{w8}
\begin{array}{c}
\dss 
\left[ \intO{ \left[ \frac{1}{2} \vr |\vu |^2 + \frac{1}{\ep^2 (\gamma - 1) } (\vr \Theta)^{\gamma} \right] } \right]_{t = 0}^{t = \tau}
+ \nu \int_{0}^{\tau} \intO{ \tn{S}(\Grad \vu) : \Grad \vu } \ \dt 
\\[20pt] 
\dss\leq \int_{0}^{\tau} \intO{ \frac{1}{\ep^2} \vr \Grad F \cdot \vu } \ \dt.
\end{array}
\eF

It was shown by Mich\' alek \cite{Mich} that the family of weak solutions to problem (\ref{w4a}), (\ref{w5a}), (\ref{w7}), and (\ref{w8}) is weakly sequentially
compact, meaning any sequence of weak solutions contains a subsequences weakly converging to another weak solution to the same problem. Novotn\' y et al. \cite{MMNP} established \emph{existence} of global-in-time weak solutions in a slightly different setting for the unknowns $[\vr, \vu, Z \equiv \vr \Theta]$, where
(\ref{w5}) is replaced by
\begin{equation}
\partial_t Z + \Div (Z \vu) = 0.
\end{equation}
For further comments concerning the existence of weak solutions satisfying (\ref{w5}) and/or (\ref{w5a}) see Section \ref{C}.

\subsection{Solutions to the target system}

As observed in the introductory part, the target system in the inviscid limit is (\ref{i5a}--\ref{i15}). In view of the global existence result by Oliver
\cite[Theorem 3]{Oli} for 2D system, we may anticipate the existence of local-in-time strong solutions $[\vc{v}, \mT, \Pi]$ also in the 3D case, specifically,
\bFormula{vspace}
\begin{array}{rcl}
\dss\vc{v} 
  & \in 
    & \dss C([0,T_{\rm max}); W^{m,2}(\Omega; \R^3))
      \\[10pt]
\dss\Pi, \mT 
  & \in 
    & \dss C([0,T_{\rm max}); W^{m,2}(\Omega)) 
\end{array}\qquad (m \geq 3)
\eF
provided
\bFormula{vstart}
\begin{array}{rcl}
\dss \vc{v}(0, \cdot) 
  & = 
    & \dss \vc{v}_0 \in W^{m,2}(\Omega; \R^3),
      \\[10pt]
\dss \mT(0, \cdot) 
  & = 
    & \dss \mT_0 \in W^{m,2}(\Omega), 
      \\[10pt]
\dss \vc{v}_0 \cdot \vc{n}|_{\partial \Omega} 
  & = 
    & \dss 0,
      \\[10pt]
\dss \Div (\tilde \vr \vc{v}_0) 
  & = 
    & 0.
\end{array}
\eF
For fixed positive viscosity coefficient $\nu > 0$, the limit problem is the Navier-Stokes system
\begin{eqnarray}
\label{w7a}
\Div (\tilde \vr \vc{v} ) 
  & = 
    & 0,
      \\
\label{w7b}
\tilde \vr \partial_t \vc{v} + \tilde \vr \vc{v} \cdot \Grad \vc{v} + \tilde \vr \Grad \Pi 
  & = 
    & \nu \Div \tn{S}(\Grad \vc{v}) - \tilde \vr \mT \Grad F,
\end{eqnarray}
where $\mT$ satisfies the transport equation (\ref{i5b}). The weak solutions of (\ref{w7a}), (\ref{w7b}), supplemented with the slip boundary condition
(\ref{i3c}) are defined in a standard way by requiring the integral identity
\bFormula{w7c}
\left[ \intO{ \tilde \vr \vc{v} } \right]_{t= \tau_1}^{t=\tau_2}
\eF
\[
= \int_{\tau_1}^{\tau_2} 
  \intO{ 
    \left[\rule{0pt}{12pt} \tilde \vr \vc{v} \cdot \partial_t \varphi 
         + \left(\rule{0pt}{12pt}\tilde \vr \vc{v} \otimes \vc{v} - \nu \tn{S}(\Grad \vc{v}) 
           \right) : \Grad \varphi 
         + \tilde \vr \mT \Grad F \cdot \varphi 
    \right] } \ \dt
\]
to holds
for any $\varphi \in \DC([0,T] \times \Ov{\Omega}; \R^3)$, $\Div(\tilde\vr\varphi)=0$, $\varphi \cdot \vc{n}|_{\partial \Omega} = 0$.

\subsection{Main results}

Having collected all the necessary preliminary material we are ready to state the two main results of this paper.
We start with the inviscid limit for the well-prepared initial data.

\bTheorem{w1}
Let $\Omega = T_{2D} \times (0,1)$ and $\gamma > \frac{3}{2}$. Let $[\vre, \vue, \The]$ be a weak solution of the primitive system (\ref{i1}--\ref{i3c})
on a time interval $(0,T)$, with
the initial data (\ref{i4}), where
\begin{equation}
\vr^{(1)}_{0,\ep}, \ \Theta^{(2)}_{0,\ep} \in L^\infty(\Omega), 
\qquad 
\vu_{0,\ep} \in L^\infty(\Omega; \R^3),
\end{equation}
\begin{equation}
\| \vr^{(1)}_{0,\ep} \|_{L^\infty(\Omega)} + \| \Theta^{(2)}_{0,\ep} \|_{L^\infty(\Omega)} + \| \vu_{0,\ep} \|_{L^\infty(\Omega; \R^3)} \leq D.
\end{equation}
Suppose that the target system (\ref{i5a}--\ref{i15}) admits a smooth solution $[\vc{v}, \mT]$ on the same time interval $[0,T]$ emanating from the initial data
\begin{equation}
\vc{v}(0, \cdot) = \vc{v}_0, \qquad \mT(0, \cdot) = \mT_0.
\end{equation}

Then
\bFormula{mw1}
\sup_{t \in [0,T]} \intO{ \left[ \vre |\vue - \vc{v}|^2 + \left| \frac{\vre  - \tilde \vr}{\ep} \right|^\gamma +  \vre \left| \frac{\The - 1}{\ep^2} - \mT \right|^2  \right] }
\eF
\[
\leq c(T,D) \left[ \ep + \nu + \intO{ | \vu_{0, \ep} - \vc{v}_0 |^2 + \left| \vr^{(1)}_{0,\ep} \right|^2 + \left| {\Theta^{(2)}_{0,\ep}} -
  \mT_0 \right|^2 } \right],
\]
where the constant $c=c(T,D)$ depends on the norm of the limit solution {$[\vc{v},{\cal T}]$} and on the size $D$ of the initial data perturbation.

\eT

\bRemark{mw1}

Theorem \ref{Tw1} yields convergence in the regime $\ep \to 0$, ${\nu} \to 0$ and for the well prepared initial data for which the expression on the
right-hand side of (\ref{mw1}) vanishes in the asymptotic limit $\ep,\ {\nu} \to 0$.

\eR

Our second result concerns the case of ill-prepared data.

\bTheorem{w2}
Let $\Omega = \R^2 \times (0,1)$ be an infinite slab and ${\gamma> {3}}$. Let $[\vre, \vue, \The]$ be a weak solution of the primitive system (\ref{i1}--\ref{i3c}), (\ref{i4a}),
on a time interval $(0,T)$, with
the initial data (\ref{i4}), where
\begin{equation}
\begin{array}{rcl}
\dss \vr^{(1)}_{0,\ep} 
  & \to 
    & \dss s_0 \,, 
      \\[10pt]
\dss \Theta^{(2)}_{0,\ep} 
  & \to 
    & \dss \mT_0 \ \mbox{weakly-(*) in}\ L^1 \cap L^\infty (\Omega) \,, 
      \\[10pt]
\dss \vu_{0,\ep} 
  & \to
    & \dss \vc{v}_0 \ \mbox{weakly-(*) in} \ L^2 \cap L^\infty(\Omega; \R^3).
\end{array}
\end{equation}
Then, passing to a suitable subsequence as the case may be, we have, {for $\ep\to\infty$},
\bFormula{mw2}
\left\{
\begin{array}{c}
\sup_{t \in (0,T)} \| \vre - \tilde \vr \|_{L^2 (\Omega)} \to 0,  \\ \\ \vue \to \vc{v}
\ \mbox{weakly in} \ L^2(0,T; W^{1,2}(\Omega;\R^3))\ \mbox{and (strongly) in}\
L^2_{\rm loc} ([0,T] \times \Ov{\Omega}; \R^3) , \\ \\ \frac{ \The - 1 }{\ep^2} \to \mT \ \mbox{weakly-(*) in}\ L^\infty((0,T) \times {\Omega}),
\end{array}
\right\}
\eF
where $[\vc{v}, \mT]$ is a weak solution to the target system (\ref{i5b}), (\ref{w7a}), (\ref{w7c}), with the initial data
$[\vc{v}_0, \mT_0]$.
\eT

The rest of the paper is basically devoted to the proofs of Theorems \ref{Tw1}, \ref{Tw2}.

\section{Relative energy, uniform bounds}
\label{R}

Similarly to Masmoudi \cite{MAS7}, Jiang and Wang \cite{WanJia}, among many others, we adapt the idea of Dafermos \cite{Daf4} based on the concept of relative (modulated) energy associated to the compressible Navier-Stokes system. In the present context, the relative energy functional reads
\bFormula{R1}
\mathcal{E}_\ep \left( \vr, \Theta, \vu \Big| r, \vc{U} \right) =
\intO{  \left[ \frac{1}{2} \vr |\vu - \vc{U}|^2 + \frac{1}{\ep^2} \Big( H(\vr \Theta) - H'(r)(\vr \Theta - r) - H(r) \Big) \right]},
\eF
where we have set
\bFormula{R11}
H(Z) = \frac{1}{\gamma - 1} Z^\gamma.
\eF

\subsection{Relative energy inequality}

As shown in \cite{FeNoJi}, any \emph{finite energy} weak solution to the Navier-Stokes system satisfies the relative energy inequality
%
\bFormula{R2}
\begin{array}{c}
\dss
\left[ \mathcal{E}_\ep \left( \vr, \Theta, \vu \Big| r, \vc{U} \right) \right]_{t = 0}^{t = \tau}
+ \nu \int\limits_{0}^{\tau} 
\int\limits_{\Omega} \tn{S} (\Grad (\vu-\vc{U})) : \Grad (\vu - \vc{U}) \ \dx \, \dt
\\
\leq
\dss \int\limits_{0}^{\tau} 
\int\limits_{\Omega}  \vr \left( \partial_t \vc{U} + \vu \cdot \Grad \vc{U} \right) \cdot \left( \vc{U} - \vu \right) 
+ \nu \, \tn{S} (\Grad \vc{U}) : \Grad (\vc{U} - \vu )
\\
\dss 
+ \frac{1}{\ep^2} 
\Biggl\{
   (r - \vr \Theta) \partial_t H'(r) + \Grad H'(r) \cdot (r \vc{U} - \vr \Theta \vu )  
\\
\dss\qquad\qquad\qquad\qquad  
     - \ \Div \vc{U} \Big( (\vr \Theta)^\gamma - r^\gamma \Big) 
     - \vr \Grad F \cdot (\vc{U} - \vu)\Biggr\}   \ \dx \, \dt
\end{array}
\eF
for all (smooth) ``test functions'' functions $r$, $\vc{U}$,
\bFormula{R2a}
(r - \tilde \vr) \in \DC([0,T] \times \Ov{\Omega}), 
\qquad \vc{U} \in \DC([0,T] \times \Ov{\Omega};\R^3)
\eF
such that
\begin{equation}
r > 0, \ \vc{U} \cdot \vc{n}|_{\partial \Omega} = 0.
\end{equation}

\bRemark{RR1}

A special form of (\ref{R2}) for a particular choice of the test functions $r$, $\vc{U}$ and under different assumptions on smoothness of the solution
$[\vr, \vu, \Theta]$ was derived by several authors, e.g., Desjardins \cite{Des1}, Germain \cite{Ger}, Mellet and Vasseur~\cite{MeVa2b}.

\eR

\subsection{Uniform estimates based on the (relative) energy inequality}
\label{UB}

The choice  $r = \tilde \vr$, $\vc{U} = 0$ in the relative energy inequality (\ref{R2}) gives rise to
\[
\intO{  \left[ \frac{1}{2} \vre |\vue|^2 + \frac{1}{\ep^2} \Big( H(\vre \Theta_\ep) - H'(\tilde \vr )(\vre \Theta_\ep - \tilde \vr) - H(\tilde \vr) \Big)(\tau, \cdot) \right]} \]
\[+ \nu \int_0^\tau \intO{ \tn{S}(\Grad \vue) : \Grad \vue } \ \dt
\]
\[
\leq \intO{  \left[ \frac{1}{2} \vr_{0,\ep} |\vu_{0,\ep} |^2 + \frac{1}{\ep^2} \Big( H(\vr_{0,\ep} \Theta_{0,\ep}) - H'(\tilde \vr )(\vr_{0,\ep} \Theta_{0,\ep} - \tilde \vr) - H(\tilde \vr) \Big) \right]}
\]
\[
- \frac{1}{\ep^2} \int_{0}^{\tau} \intO{ \Grad H'(\tilde \vr) \cdot ( \vre \Theta_\ep \vue ) \Big] } \ \dt
+ \frac{1}{\ep^2} \int_{0}^{\tau} \intO{ \vre \Grad F \cdot \vue } \ \dt.
\]
In view of (\ref{R11}) and
\bFormula{tr}
\Grad H'(\tilde \vr) = \Grad F
\eF
we may infer
\bFormula{B1}
\intO{  \left[ \frac{1}{2} \vre |\vue|^2 + \frac{1}{\ep^2} \Big( H(\vre \Theta_\ep) - H'(\tilde \vr )(\vre \Theta_\ep - \tilde \vr) - H(\tilde \vr) \Big)(\tau, \cdot) \right]} 
\eF
\[+ \nu \int_0^\tau \intO{ \tn{S}(\Grad \vue) : \Grad \vue } \ \dt\]
\[
\leq \intO{  \left[ \frac{1}{2} \vr_{0,\ep} |\vu_{0,\ep} |^2 + \frac{1}{\ep^2} \Big( H(\vr_{0,\ep} \Theta_{0,\ep}) - H'(\tilde \vr )(\vr_{0,\ep} \Theta_{0,\ep} - \tilde \vr) - H(\tilde \vr) \Big) \right]}
\]
\[
\frac{1}{\ep^2} \int_{0}^{\tau} \intO{ \vre \Grad F \cdot \vue (1 - \Theta_\ep) } \dt
\]
\[
\leq
\intO{  \left[ \frac{1}{2} \vr_{0,\ep} |\vu_{0,\ep} |^2 + \frac{1}{\ep^2} \Big( H(\vr_{0,\ep} \Theta_{0,\ep}) - H'(\tilde \vr )(\vr_{0,\ep} \Theta_{0,\ep} - \tilde \vr) - H(\tilde \vr) \Big) \right]}
\]
\[
+ \| \Grad F \|_{L^\infty(\Omega)} \int_0^\tau \intO{ \frac{1}{2} \left[ \vre |\vue|^2 + \vre \left( \frac{\Theta_\ep - 1}{\ep^2} \right)^2 } \right] \dt.
\]

Finally, taking $\varphi = 1$, $G = \left( \frac{\Theta_\ep - 1 }{\ep^2} \right)^2$ in (\ref{w5}) we obtain
\bFormula{B2}
\intO{ \vre \left( \frac{\Theta_\ep - 1 }{\ep^2} \right)^2 (\tau, \cdot) } \leq \intO{ \vr_{0,\ep} \left( \Theta^{(2)}_{0,\ep} \right)^2  }.
\eF

Relations (\ref{B1}), (\ref{B2}) yield uniform bounds on the family $\{ \vre, \Theta_\ep, \vue \}_{\ep > 0}$ provided the initial data
are bounded as in the hypotheses of Theorems \ref{Tw1}, \ref{Tw2}, specifically,
are chosen in such a way that
\bFormula{B3}
\| \vu_{0,\ep} \|_{L^1 \cap L^\infty(\Omega; \R^3)} \leq c,
\eF
\bFormula{B4}
\vr_{0,\ep} = \tilde \vr + \ep \vr^{(1)}_\ep , \ \| \vr^{(1)}_{0,\ep} \|_{L^1 \cap L^\infty(\Omega)} \leq c,
\eF
\bFormula{B5}
\Theta_{0,\ep} = 1 + \ep^2 \Theta^{(2)}_{0,\ep} ,\ \| \Theta^{(2)}_{0,\ep} \|_{L^1 \cap L^\infty(\Omega)} \leq c.
\eF
Now, we deduce from (\ref{B1}), (\ref{B2}) the following uniform bounds independent of $\ep$, $\nu$:
\bFormula{B6}
{\rm ess} \sup_{t \in [0,T]} \| \sqrt{\vre} \vue \|_{L^2(\Omega; \R^3)} \leq c,
\eF
\bFormula{B7}
\sqrt{\nu} \left\| \Grad \vue + \Grad^t \vue - \frac{2}{3} \Div \vue \tn{I} \right\|_{L^2((0,T) \times \Omega)} \leq c,
\eF
\bFormula{B8}
{\rm ess} \sup_{t \in (0,T)} \left\| \sqrt{\vre} \left( \frac{ \Theta_\ep - 1 }{\ep^2} \right) \right\|_{L^2(\Omega)} \leq c,
\eF
\bFormula{B9}
{\rm ess} \sup_{t \in (0,T)} \int_{ \tilde \vr/ 2 \leq \vre \Theta_\ep \leq 2 \tilde \vr} \left| \frac{ \vre \Theta_\ep - \tilde \vr }{\ep} \right|^2 \ \dx \leq c,
\eF
\bFormula{B10}
{\rm ess} \sup_{t \in (0,T)} \int_{ \vre \Theta_\ep  < \tilde \vr/2} \left(1 +  \left| { \vre \Theta_\ep }  \right|^\gamma \right) \ \dx +
{\rm ess} \sup_{t \in (0,T)} \int_{ \vre \Theta_\ep >  2 \tilde \vr} \left(1 +  \left| { \vre \Theta_\ep }  \right|^\gamma \right) \ \dx
\leq \ep^2 c.
\eF

As a matter of fact, relation (\ref{w5}) may be used repeatedly to deduce
\bFormula{B11}
{\rm ess} \sup_{t \in (0,T)} \left\|  \frac{ \Theta_\ep - 1 }{\ep^2}  \right\|_{L^\infty(\{ \vre(t, \cdot) > 0 \})} \leq c
\eF
as well as
\bFormula{B12}
{\rm ess} \sup_{t \in (0,T)} \intO{ \vre \left| \frac{\Theta_\ep - 1}{\ep^2} \right| } \leq c.
\eF

\bRemark{RR2}

As the function $\Theta_\ep$ appears in (\ref{i1}--\ref{i3}) only multiplied by $\vre$, we may assume, without loss of generality, that
$\Theta_\ep = 1$ on the (hypothetical) vacuum set $\vre = 0$. Keeping this convention in mind, we observe that (\ref{B11}) holds in the whole physical space $\Omega$.

\eR

Now, we have
\[
\int_{ \tilde \vr/2 \leq \vre \Theta_\ep \leq 2 \tilde \vr } \left| \frac{\vre - \tilde \vr}{\ep} \right|^2 \ \dx \leq
\]
\[
\int_{ \tilde \vr/2 \leq \vre \Theta_\ep \leq 2 \tilde \vr } \left| \frac{\vre(\The - 1)}{\ep} \right|^2 \ \dx
+ \int_{ \tilde \vr/2 \leq \vre \Theta_\ep \leq 2 \tilde \vr } \left| \frac{\vre \The - \tilde \vr}{\ep} \right|^2 \ \dx,
\]
therefore (\ref{B8}), (\ref{B9}) give rise to
\bFormula{B13}
{\rm ess}\sup_{t \in (0,T)}\int_{ \tilde \vr/2 \leq \vre \Theta_\ep \leq 2 \tilde \vr } \left| \frac{\vre - \tilde \vr}{\ep} \right|^2 \ \dx \leq c.
\eF
Similarly, the bounds (\ref{B10}), (\ref{B11}) yield
\bFormula{B14}
{\rm ess}\sup_{t \in (0,T)} \int_{ \vre \The < \tilde \vr/2 } |\vre|^\gamma \ \dx +
{\rm ess}\sup_{t \in (0,T)} \int_{ \vre \The > 2 \tilde \vr } |\vre|^\gamma \ \dx \leq \ep^2 c.
\eF

\
At this stage it is convenient to introduce the following notation:
for each measurable function $f(t,x)$ we write $f=[f]_{{\rm ess}}+[f]_{{\rm res}}$, where
\bFormula{ess-res}
[f]_{{\rm ess}}=\chi(\vr_\ep\Theta_\ep)f,\qquad [f]_{{\rm res}}=(1-\chi(\vre\Theta_\ep))f,
\eF
\[\chi=\left\{\begin{array}{l}1\quad \mbox{for}\  \tilde\vr/2\leq\vre\Theta_\ep\leq 2\tilde\vr\\
0\quad \mbox{otherwise.}
\end{array}\right.
\]
Now, we have
\bFormula{B15a}
\intO{ \tilde \vr \left| \frac{\The - 1}{\ep^2} \right| } \leq \intO{ \left|\tilde \vr - \frac{\vre\Theta_\ep}{2} \right| \left| \frac{ \The - 1}{\ep^2} \right| }
+ \frac{1}{2}\intO{ \vre\Theta_\ep  \left| \frac{ \The - 1}{\ep^2} \right| } ;
\eF
where the last integral is bounded due to (\ref{B11}) and (\ref{B12}). As for the first integral, we may decompose the integrant into its {\it essential} and {\it residual} part, where
the essential part can be estimated as follows
\bFormula{B15b}
\begin{array}{c}
\dss
 \intO{\left[\left|\tilde \vr - \frac{\vre \Theta_\ep}{2}\right| \left| \frac{ \The - 1}{\ep^2} \right|\right]_{{\rm ess}}}
 \leq \frac{3}{4} \intO{\left[\tilde \vr  \left| \frac{ \The - 1}{\ep^2} \right|\right]_{{\rm ess}}}\\
 \dss
 \hspace{5cm}\leq \frac{3}{4}\intO{\tilde \vr  \left| \frac{ \The - 1}{\ep^2} \right|},
 \end{array}
\eF
while for the residual part we have
\bFormula{B15c}
\begin{array}{c}
\dss
 \intO{\left[ \left|\tilde \vr - \frac{\vre \Theta_\ep}{2}\right| \left| \frac{ \The - 1}{\ep^2} \right|\right]_{{\rm res}}}\\
 \dss
\hspace{4cm} \leq   \left\| \frac{ \The - 1}{\ep^2} \right\|_{L^\infty(\Omega)}\intO{\left[\left|\tilde \vr - \frac{\vre \Theta_\ep}{2}\right|\right]_{{\rm res}}}
 \leq c
 \end{array}
\eF
which follows from estimates (\ref{B10}) and (\ref{B11}) and the fact that $\|\tilde\vr\|_{L^\infty(\Omega)}\leq c.$
Putting together estimates (\ref{B15a}-\ref{B15c}) and recalling that $\tilde\rho>0$, we deduce
\bFormula{B15}
{\rm ess}\sup_{t \in (0,T)} \left\| \frac{ \The - 1 }{\ep^2} \right\|_{L^1(\Omega)} \leq c.
\eF


\section{The inviscid limit}
\label{A}

In this section, we prove Theorem \ref{Tw1}. To this end, we take $\vc{U} = \vc{v}$, $r = \tilde \vr$ in as test functions in the relative
energy inequality (\ref{R2}), where $\vc{v}$ is the (smooth) solution of the target system (\ref{i5a}--\ref{i15}). After a straightforward manipulation (see \cite{FeLuNo} for more details) we
obtain
\bFormula{A4}
\left[ \mathcal{E}_\ep \left( \vre, \Theta_\ep, \vue \Big| \tilde \vr , \vc{v} \right) \right]_{t = 0}^{t = \tau}
+ \frac{\nu}{2} \int_{\tau_1}^{\tau_2} \intO{ \Big( \tn{S} (\Grad \vue) - \tn{S} (\Grad \vc{v}) \Big) : \Big( \Grad \vue - \Grad \vc{v} \Big) } \ \dt
\eF
\[
\leq
\int_{0}^{\tau} \intO{  \vre \left( \partial_t \vc{v} + \vue \cdot \Grad \vc{v} \right) \cdot \left( \vc{v} - \vue \right) } \ \dt
+ \nu \int_{0}^{\tau} \intO{ |\Grad \vc{v} |^2 } \ \dt
\]
\[
+ \frac{1}{\ep^2} \int_{0}^{\tau} \intO{ \Big[ (\tilde \vr  - \vre \Theta_\ep) \partial_t H'(\tilde \vr ) + \Grad H'(\tilde \vr ) \cdot (\tilde \vr \vc{v} - \vre \Theta_\ep \vue ) \Big] } \ \dt
\]
\[
- \frac{1}{\ep^2} \int_{0}^{\tau} \intO{ \Div \vc{v} \Big( (\vre \Theta_\ep)^\gamma - \tilde \vr^\gamma \Big) } \ \dt
- \frac{1}{\ep^2} \int_{0}^{\tau} \intO{ \vre \Grad F \cdot (\vc{v} - \vue) } \ \dt
\]
\[
\leq \| \Grad \vc{v} \|_{L^\infty((0,T) \times \Omega; \R^{3 \times 3} )} \int_0^T \intO{ \vre |\vue - \vc{v} |^2 } +
\nu \int_{0}^{\tau} \intO{ |\Grad \vc{v} |^2 } \ \dt
\]
\[
\int_{0}^{\tau} \intO{  \vre \left( - \mT \Grad F - \Grad \Pi \right) \cdot \left( \vc{v} - \vue \right) } \ \dt
\]
\[
+ \frac{1}{\ep^2} \int_{0}^{\tau} \intO{  \Grad H'(\tilde \vr) \cdot (\tilde \vr  \vc{v} - \vre \Theta_\ep \vue ) } \ \dt
\]
\[
- \frac{1}{\ep^2} \int_{0}^{\tau} \intO{ \Div \vc{v} \Big( (\vre \Theta_\ep)^\gamma - (\tilde \vr \mT)^\gamma \Big) } \ \dt
- \frac{1}{\ep^2} \int_{0}^{\tau} \intO{ \vre \Grad F \cdot (\vc{v} - \vue) } \ \dt
\]
\[
= \| \Grad \vc{v} \|_{L^\infty((0,T) \times \Omega; \R^{3 \times 3} )} \int_0^T \intO{ \vre |\vue - \vc{v} |^2 } +
\nu \int_{0}^{\tau} \intO{ |\Grad \vc{v} |^2 } \ \dt
\]
\[
\int_{0}^{\tau} \intO{  \vre  \left( - \mT - \frac{1}{\ep^2} \right) \Grad F \cdot \left( \vc{v} - \vue \right) } \ \dt
+ \int_0^\tau \intO{ \vre \Grad \Pi \cdot (\vue - \vc{v} ) } \ \dt
\]
\[
+ \frac{1}{\ep^2} \int_{0}^{\tau} \intO{  \Grad H'(\tilde \vr) \cdot (\tilde \vr  \vc{v} - \vre \Theta_\ep \vue ) } \ \dt
- \frac{1}{\ep^2} \int_{0}^{\tau} \intO{ \Div \vc{v} \Big( (\vre \Theta_\ep)^\gamma - \tilde \vr^\gamma \Big) } \ \dt.
\]

Next, we use (\ref{tr}) to rewrite the penultimate term as follows
\[
\frac{1}{\ep^2} \int_{0}^{\tau} \intO{  \Grad H'(\tilde \vr) \cdot (\tilde \vr  \vc{v} - \vre \Theta_\ep \vue ) } \ \dt
\]
\[
=\frac{1}{\ep^2} \int_{0}^{\tau} \intO{  (\tilde \vr  - \vre \Theta_\ep  ) \Grad H'(\tilde \vr) \cdot \vc{v}  } \ \dt
+ \frac{1}{\ep^2} \int_{0}^{\tau} \intO{ \vre \Theta_\ep   \Grad H'(\tilde \vr) \cdot ( \vc{v} - \vue)   } \ \dt
\]
\[
= \frac{1}{\ep^2} \int_{0}^{\tau} \intO{  (\tilde \vr  - \vre \Theta_\ep  ) \Grad H'(\tilde \vr) \cdot \vc{v}  } \ \dt
+ \frac{1}{\ep^2} \int_{0}^{\tau} \intO{ \vre \Theta_\ep   \Grad F \cdot ( \vc{v} - \vue)   } \ \dt.
\]

Furthermore, in accordance with the anelastic constraint (\ref{i5a}),
\[
\Grad \tilde \vr \cdot \vc{v} = - \tilde \vr \Div \vc{v};
\]
whence
\[
\frac{1}{\ep^2} \int_{0}^{\tau} \intO{  (\tilde \vr  - \vre \Theta_\ep  ) \Grad H'(\tilde \vr) \cdot \vc{v}  } \ \dt =
\frac{1}{\ep^2} \int_{0}^{\tau} \intO{  (\vre \Theta_\ep -\tilde \vr  ) \gamma \tilde \vr^{\gamma - 1} \Div \vc{v}  } \ \dt.
\]
Thus, revisiting (\ref{A4}), we deduce
\bFormula{A5}
\left[ \mathcal{E}_\ep \left( \vre, \Theta_\ep, \vue \Big| \tilde \vr , \vc{v} \right) \right]_{t = 0}^{t = \tau}
+ \frac{\nu}{2} \int_{0}^{\tau} \intO{ \Big( \tn{S} (\Grad \vue) - \tn{S} (\Grad \vc{v}) \Big) : \Big( \Grad \vue - \Grad \vc{v} \Big) } \ \dt
\eF
\[
\leq
\| \Grad \vc{v} \|_{L^\infty((0,T) \times \Omega; \R^{3 \times 3} )} \int_0^\tau \intO{ \mathcal{E}_\ep \left( \vre, \Theta_\ep, \vue \Big| \tilde \vr , \vc{v} \right) } \dt +
\nu \int_{0}^{\tau} \intO{ |\Grad \vc{v} |^2 } \ \dt
\]
\[
+\int_{0}^{\tau} \intO{  \vre  \left( \frac{\Theta_\ep - 1}{\ep^2} - \mT \right) \Grad F \cdot \left( \vc{v} - \vue \right) } \ \dt
+ \int_0^\tau \intO{ \vre \Grad \Pi \cdot (\vue - \vc{v} ) } \ \dt.
\]
The last integral on the right-hand side can be handled as follows:
\[
\int_0^\tau \intO{ \vre \Grad \Pi \cdot (\vue - \vc{v} ) } \ \dt
\]
\[
= - \int_0^\tau \intO{ \vre \partial_t \Pi } +
\left[ \intO{ \vre \Pi } \right]_{t = 0}^{t = \tau} + \int_0^\tau \intO{ (\tilde \vr - \vre) \Grad \Pi \cdot \vc{v} };
\]
whence
\bFormula{A6}
\left| \int_0^\tau \intO{ \vre \Grad \Pi \cdot (\vue - \vc{v} ) } \ \dt \right| \leq c \ {\rm ess} \sup_{t \in [0,T]} \| \vre - \tilde \vr \|_{L^1(\Omega)}.
\eF
Going back to (\ref{A5}) we obtain
\bFormula{A7}
\left[ \mathcal{E}_\ep \left( \vre, \Theta_\ep, \vue \Big| \tilde \vr , \vc{v} \right) \right]_{t = 0}^{t = \tau}
\eF
\[
\leq
c \left( \| \Grad \vc{v} \|_{L^\infty((0,T) \times \Omega; \R^{3 \times 3} )}  + \| \Grad F \|_{L^\infty((0,T) \times \Omega; \R^{3} )} \right)
\int_0^\tau \intO{ \mathcal{E}_\ep \left( \vre, \Theta_\ep, \vue \Big| \tilde \vr , \vc{v} \right) } \dt
\]
\[
+
\nu \int_{0}^{\tau} \intO{ |\Grad \vc{v} |^2 } \ \dt
\]
\[
+c \| \Grad F \|_{L^\infty((0,T) \times \Omega; \R^{3} )}
\int_{0}^{\tau} \intO{  \vre  \left( \frac{\Theta_\ep - 1}{\ep^2} - \mT \right)^2  } \ \dt\]
\[+ c \| \Pi \|_{W^{1,\infty}((0,T) \times \Omega)} \ {\rm ess} \sup_{t \in [0,T]} \| \vre - \tilde \vr \|_{L^1(\Omega)}.
\]
Finally, to bound the penultimate term we show:
\bLemma{w1}
Suppose $\vr$, $\Theta$ is a weak solution of (\ref{i1}), (\ref{i3}) in the sense specified in (\ref{w5}). Let
\bFormula{w6-}
\partial_t \mT + \vc{v} \cdot \Grad \mT = 0, \ \mT \in C^1, \vc{v} \in C^1.
\eF

Then
\bFormula{w6}
\left[ \intO{ \frac{1}{2} \vr | G(\Theta) - \mT |^2 } \right]_{t = \tau_1}^{t = \tau_2}
= \int_{\tau_1}^{\tau_2} \intO{  \vr \Big( G(\Theta) - \mT \Big) \Big(\vu -  \vc{v} \Big) \cdot \Grad \mT } \ \dt.
\eF
\eL

\bProof

We write
\[
\left[ \intO{ \frac{1}{2} \vr | G(\Theta) - \mT |^2 } \right]_{t = \tau_1}^{t = \tau_2}
\]
\[
= \left[ \intO{ \frac{1}{2} \vr G(\Theta)^2 } \right]_{t = \tau_1}^{t = \tau_2} + \left[ \intO{ \vr G(\Theta) \mT } \right]_{t = \tau_1}^{t = \tau_2}
+ \left[ \intO{ \frac{1}{2} \vr \mT^2 } \right]_{t = \tau_1}^{t = \tau_2},
\]
where, in accordance with (\ref{w5}),
\[
\left[ \intO{ \frac{1}{2} \vr G(\Theta)^2 } \right]_{t = \tau_1}^{t = \tau_2} = 0.
\]

Next, using the weak formulation (\ref{w5}) with the test functions $\varphi  =\mT$ and again with $G \equiv 1$, $\varphi = \frac{1}{2} \mT^2$, we obtain
\[
\left[ \intO{ \vr G(\Theta) \mT } \right]_{t = \tau_1}^{t = \tau_2}
+ \left[ \intO{ \frac{1}{2} \vr \mT^2 } \right]_{t = \tau_1}^{t = \tau_2}
\]
\[
=\int_{\tau_1}^{\tau_2} \intO{ \Big[ \vr G(\Theta) \partial_t \mT + \vr G(\Theta) \vu \cdot \Grad \mT + \vr \mT \partial_t \mT + \vr  \vu \cdot \mT \Grad \mT \Big] } \ \dt;
\]
whence, by virtue of (\ref{w6-}),
\[
\left[ \intO{ \vr G(\Theta) \mT } \right]_{t = \tau_1}^{t = \tau_2}
+ \left[ \intO{ \frac{1}{2} \vr \mT^2 } \right]_{t = \tau_1}^{t = \tau_2}
\]
\[
=\int_{\tau_1}^{\tau_2} \intO{ \Big[ - \vr G(\Theta) \vc{v} \cdot \Grad \mT + \vr G(\Theta) \vu \cdot \Grad \mT - \vr \mT \vc{v} \cdot \Grad \mT + \vr  \vu \cdot \mT \Grad \mT \Big] } \ \dt
\]
\[
= \int_{\tau_1}^{\tau_2} \intO{  \vr \Big( G(\Theta) - \mT \Big) \Big(\vu -  \vc{v} \Big) \cdot \Grad \mT } \ \dt.
\]

\qed

Using (\ref{w6}) and the uniform bounds established in (\ref{B6}--\ref{B11}), we may now apply Gronwall's lemma to (\ref{A7}), to obtain convergence towards the limit system.
In addition, the rate of convergence can be computed explicitly in terms of $\ep$, $\nu$ as indicated in (\ref{mw1}). Theorem \ref{Tw1} has been proved.

\section{The ill-prepared initial data}
\label{I}

With $\nu > 0$ fixed, our ultimate goal is to preform the singular limit $\ep \to 0$ for the ill-prepared data as claimed in Theorem \ref{Tw2}. This
is definitely more delicate than the previous part as the motion is polluted by rapidly oscillating acoustic waves. Fortunately, the unbounded physical spaces
gives rise to dispersive effects, in particular, the acoustic part of the velocity vanishes in the asymptotic limit, at least on compact subsets of $\Omega$.

\subsection{Weak convergence}

To begin, a generalized version of Korn's inequality (see e.g. \cite[Section 10.9]{FeNo6}) can be used to deduce from (\ref{B7}), (\ref{B13}), and (\ref{B14}) the
bound
\bFormula{wc1}
\int_0^T \intO{ |\Grad \vue |^2 + |\vue |^2 } \leq c,
\eF
in particular, we may assume that
\[
\vue \to \vc{v} \ \mbox{weakly in}\ L^2(0,T; W^{1,2}(\Omega;\R^3)),
\]
passing to a subsequence if necessary. Similarly, it is easy to check that
\[
{\rm ess} \sup_{t \in (0,T)} \| \vre - \tilde \vr \|_{L^2(\Omega; \R^3)} \to 0
\]
keeping in mind that $\gamma \geq 3$ in Theorem \ref{Tw2}.

Finally, it is easy to observe that the continuity equation (\ref{i1}) reduces to the anelastic constraint (\ref{i5a}) in the asymptotic limit $\ep \to 0$.

\subsection{Momentum balance}

We write the momentum equation (\ref{i2}) in the form
\bFormula{AA2-}
\begin{array}{c}
\dss
\partial_t (\vre \vue) + \frac{1}{\ep^2} \left( \Grad (\vre \Theta_\ep )^\gamma - \vre \Theta_\ep \Grad F \right)\\
\dss =
- \Div (\vre \vue \otimes \vue ) {+ \nu\Div \tn{S}(\Grad \vue)} + \vre \frac{1 - \The}{\ep^2} \Grad F,
\end{array}
\eF
where, furthermore,
\[
\Grad (\vre \Theta_\ep )^\gamma - \vre \Theta_\ep \Grad F
\]
\[
= \Grad \left[ (\vre \Theta_\ep)^\gamma  - \gamma {\tilde \vr}^{\gamma - 1}(
\vre \Theta_\ep - \tilde \vr) - {\tilde \vr}^\gamma \right] \]
\[- \gamma \vre \Theta_\ep {\tilde \vr}^{\gamma - 2} \Grad \tilde \vr
+  \gamma {\tilde \vr}^{\gamma - 1} \Grad \tilde \vr + \Grad \left[ \gamma {\tilde \vr}^{\gamma - 1}(
\vre \Theta_\ep - \tilde \vr) \right]
\]
\[
= \Grad \left[ (\vre \Theta_\ep)^\gamma  - \gamma {\tilde \vr}^{\gamma - 1}(
\vre \Theta_\ep - \tilde \vr) - {\tilde \vr}^\gamma \right]- ( \vre \Theta_\ep - \tilde \vr) \gamma {\tilde \vr}^{\gamma - 2} \Grad \tilde \vr
 + \Grad \left[ \gamma {\tilde \vr}^{\gamma - 1}(
\vre \Theta_\ep - \tilde \vr) \right]
\]
\[
= \Grad \left[ (\vre \Theta_\ep)^\gamma  - \gamma {\tilde \vr}^{\gamma - 1}(
\vre \Theta_\ep - \tilde \vr) - {\tilde \vr}^\gamma \right] + \tilde \vr \Grad \left[ \gamma {\tilde \vr}^{\gamma - 2} \left( \vre \Theta_\ep - \tilde \vr \right)     \right].
\]

Consequently, a suitable weak formulation of (\ref{AA2-}) reads
\bFormula{AA2+}
\int_0^T \intO{ \left[
\vre \vue \cdot \partial_t \varphi  + \vre \vue \otimes \vue : \Grad \varphi \right] } \ \dt
\eF
\[
\int_0^T \intO{ \frac{ (\vre \Theta_\ep)^\gamma  - \gamma {\tilde \vr}^{\gamma - 1}(
\vre \Theta_\ep - \tilde \vr) - {\tilde \vr}^\gamma }{\ep^2} \Div \varphi} \ \dt\] 
\[+\int_0^T \intO{
\frac{1}{\ep^2} \left( \gamma {\tilde \vr}^{\gamma - 2} \left( \vre \Theta_\ep - \tilde \vr \right) \right) \Div( \tilde \vr \varphi)
} \ \dt\]
\[
=
\int_0^T \intO{ \left[ {\nu\tn{S}(\Grad \vue)}: \Grad \varphi  - \vre \frac{1 - \The}{\ep^2} \Grad F \cdot \varphi \right] } \ \dt
- \intO{ \vr_{0,\ep} \vu_{0,\ep} \cdot \varphi (0, \cdot) }
\]
for any $\varphi \in \DC([0,T) \times \Ov{\Omega}$, $\varphi \cdot \vc{n}|_{\partial \Omega} = 0$.

Next, in agreement with the uniform bounds
(\ref{B11}), (\ref{B15}), we may assume that
\[
{\frac{ \The - 1}{\ep^2}} \to \mT \ \mbox{weakly-(*) in}\ L^\infty(0,T; L^1 \cap L^\infty(\Omega)).
\]

Consequently, it is possible to perform the limit for $\ep \to 0$ in both (\ref{AA2+}), for $\varphi$ satisfying $\Div (\tilde \vr \varphi) = 0$, and
(\ref{i3}) to obtain
\bFormula{AA2++}
\int_0^T \intO{ \left[
\tilde \vr \vc{v} \cdot \partial_t \varphi  + \Ov{ \vr \vu \otimes \vu} : \Grad \varphi \right] } \ \dt
\eF
\[
\lim_{\ep \to 0} \int_0^T \intO{ \frac{ (\vre \Theta_\ep)^\gamma  - \gamma {\tilde \vr}^{\gamma - 1}(
\vre \Theta_\ep - \tilde \vr) - {\tilde \vr}^\gamma }{\ep^2} \Div \varphi } \ \dt
\]
\[
=
\int_0^T \intO{ \left[ {\nu\tn{S}(\Grad \vc{v})}: \Grad \varphi  - \tilde \vr  \mT \Grad F \cdot \varphi \right] } \ \dt
- \intO{ \tilde \vr \vc{v}_{0} \cdot \varphi (0, \cdot) }
\]
for any $\varphi \in \DC([0,T) \times \Ov{\Omega}$, $\varphi \cdot \vc{n}|_{\partial \Omega} = 0$, $\Div (\tilde \vr \varphi ) = 0$,
and
\bFormula{i3ca}
\int_0^T \intO{ \left[ \tilde \vr \mT \partial_t \varphi  + \tilde \vr \mT \vc{v} \cdot \Grad \varphi \right] } \ \dt = - \intO{ \tilde \vr \mT_0 \varphi (0, \cdot) }
\eF
for any $\varphi \in \DC([0,T) \times \Ov{\Omega})$, where
$\Ov{\vr \vu \otimes \vu}$ denotes a weak limit of the sequence $\{ \vre \vue \otimes \vue \}_{\ep > 0}$.

\bRemark{Ri1}

Note that we have used Lions-Aubin lemma to obtain {(\ref{i3ca})}.

\eR

Thus, in order to complete the proof of Theorem \ref{Tw2}, it remains to show:

\begin{itemize}
\item
\bFormula{final1}
\lim_{\ep \to 0} \int_0^T \intO{ \frac{ (\vre \Theta_\ep)^\gamma  - \gamma {\tilde \vr}^{\gamma - 1}(
\vre \Theta_\ep - \tilde \vr) - {\tilde \vr}^\gamma }{\ep^2} \Div \varphi } \ \dt \to 0,
\eF
\item
\bFormula{final2}
\Ov{ \vr \vu \otimes \vu} = \tilde \vr \vc{v} \otimes \vc{v}.
\eF

\end{itemize}
This will be done in the next section by means of careful analysis of propagation of acoustic waves.

\subsection{Acoustic equation}

We start by introducing a generalized Helmholtz decomposition
\bFormula{helmholtz}
\vc{w} = \mathcal{P}_{\tilde \vr} (\vc{w}) + \tilde \vr \mathcal{Q}_{\tilde \vr} (\vc{w}), \qquad \mathcal{Q}_{\tilde \vr}(\vc{w}) = \Grad \Psi,
\eF
where $\Psi$ is the unique solution of the Neumann problem
\[
\Div (\tilde \vr \Grad \Psi)  = \Div \vc{w}, \ (\vc{w} - \tilde \vr \Grad \Psi) \cdot \vc{n} |_{\partial \Omega} = 0, \ \Psi \to 0 \ \mbox{as} \ |x| \to \infty.
\]

Accordingly, we write
\[
\vre \vue = \mathcal{P}_{\tilde \vr} (\vre \vue) + \tilde \vr \mathcal{Q}_{\tilde \vr} (\vre \vue).
\]
Note that, by virtue of (\ref{B6}), (\ref{B14}), (\ref{wc1}), and the hypothesis {$\gamma > 3$}, we may write
\[
\vre \vue = [\vre \vue]_{\rm ess} + [\vre \vue]_{\rm res},
\]
cf. (\ref{ess-res}). We have
\[
[\vre \vue]_{\rm ess} \to \tilde \vr \vc{v} \ \mbox{weakly-(*) in}\ L^\infty(0,T; L^2(\Omega;\R^3)),
\]
while
\[
[\vre \vue]_{\rm res} \to 0 \ \mbox{in}\ L^2((0,T) \times \Omega; \R^3).
\]

Seeing that
\[
\intO{ \vc{w} \cdot \varphi } = \intO{ \mathcal{P}_{\tilde \vr}(\vc{w}) \cdot \varphi } \ \mbox{whenever}\ \Div (\tilde \vr \vc{w}) = 0,
\]
we may use the momentum equation (\ref{AA2+}), together with the bounds established in Section \ref{UB}, (\ref{wc1}), and the standard Lions-Aubin argument
to conclude that
\bFormula{Aa1}
\mathcal{P}_{\tilde \vr} (\vre \vue) \to \mathcal{P}_{\tilde \vr}(\tilde \vr \vc{v}) = \tilde \vr \vc{v} \ \mbox{(strongly) in}\ L^2_{\rm loc} ((0,T) \times \Omega; \R^3).
\eF

To derive an equation for the acoustic potential we introduce notation
\[
S_\ep = \frac{ \vre \Theta_\ep - \tilde \vr }{\ep \tilde \vr}, \qquad \Grad \Phi_\ep = \mathcal{Q}_{\tilde \vr}(\vre \vue).
\]
First, we rewrite equation (\ref{i3}) in the form
\bFormula{AA3}
\ep \partial_t S_\ep + \frac{1}{\tilde \vr} \Div ( \tilde \vr \Grad \Phi_\ep) = \ep \frac{1}{\tilde \vr} \Div \left( \vre
\left( \frac{1 - \Theta_\ep}{\ep} \right) \vue \right),
\eF
whereas (\ref{AA2+}) reads
\bFormula{AA4}
\ep \partial_t (\vre \vue) + \tilde \vr \Grad \left[ \gamma {\tilde \vr}^{\gamma - 1} S_\ep \right]
\eF
\[
= - \ep \Grad \left[ \frac{(\vre \Theta_\ep)^\gamma  - \gamma {\tilde \vr}^{\gamma - 1}(
\vre \Theta_\ep - \tilde \vr) - {\tilde \vr}^\gamma}{\ep^2} \right] - \ep \Div (\vre \vue \otimes \vue ) \]
\[{+\ep \nu\Div \tn{S}(\Grad \vue) }+
\ep \vre \frac{1 - \The}{\ep^2} \Grad F.
\]

Our goal is to apply the generalized Helmholtz projection $\mathcal{Q}_{\tilde \vr}$ to (\ref{AA4}). This step may be performed in
the weak sense by taking a test function
\[
\Grad \phi \in \DC((0,T) \times \Ov{\Omega}; \R^3), \ \Grad \phi \cdot \vc{n}|_{\partial \Omega} = 0
\]
in the weak formulation of (\ref{AA4}):
\bFormula{AA7}
\int_0^T \intO{ \left[ \ep \tilde \vr \Grad \Phi_\ep \cdot \partial_t \Grad \phi + \left( \gamma {\tilde \vr}^{\gamma - 1} S_\ep \right) \Div (\tilde \vr \Grad \phi ) \right]
} \ \dt
\eF
\[
= - \ep \int_0^T \intO{ \left[ \frac{(\vre \Theta_\ep)^\gamma  - \gamma {\tilde \vr}^{\gamma - 1}(
\vre \Theta_\ep - \tilde \vr) - {\tilde \vr}^\gamma}{\ep^2} \right] \Delta \phi } \ \dt
\]
\[
- \ep \int_0^T \intO{ \left[ (\vre \vue \otimes \vue ): \Grad^2 \phi - \nu\tn{S}(\Grad \vue) : \Grad^2 \phi
+ \vre \frac{1 - \The}{\ep^2} \Grad F \cdot \Grad \phi \right] } \ \dt.
\]
Now, for $\varphi \in \DC((0,T) \times \Ov{\Omega}; \R^3$, $\varphi \cdot \vc{n}|_{\partial \Omega} = 0$, we may take $\Grad \phi = Q_{\tilde \vr}(\varphi)$ as
a test function in (\ref{AA7}) to obtain:
\bFormula{AA8}
\int_0^T \intO{ \left[ - \ep \Phi_\ep \partial_t \Div \varphi + \left( \gamma {\tilde \vr}^{\gamma - 1} S_\ep \right) \Div \varphi \right]
} \ \dt
\eF
\[
= - \ep \int_0^T \intO{ \left[ \frac{(\vre \Theta_\ep)^\gamma  - \gamma {\tilde \vr}^{\gamma - 1}(
\vre \Theta_\ep - \tilde \vr) - {\tilde \vr}^\gamma}{\ep^2} \right] \Div \mathcal{Q}_{\tilde \varrho} (\varphi) } \ \dt
\]
\[
- \ep \int_0^T \intO{ \left[ (\vre \vue \otimes \vue ): \Grad \mathcal{Q}_{\tilde \vr} (\varphi) - \nu\tn{S}(\Grad \vue) : \Grad \mathcal{Q}_{\tilde \vr} (\varphi) \right] } \ \dt
\]
\[-\ep\int_0^T \intO{\vre \frac{1 - \The}{\ep^2} \Grad F \cdot \mathcal{Q}_{\tilde \vr} (\varphi)}\ \dt\]

We conclude that the system (\ref{AA3}), (\ref{AA8}) can be written as a variant of \emph{Lighthill's acoustic analogy}:
\bFormula{AA9}
\ep \partial_t S_\ep + \frac{1}{\tilde \vr} \Div ( \tilde \vr \Grad \Phi_\ep) = \ep \mathcal{G}^1_\ep,
\eF
\bFormula{AA10}
\ep \partial_t \Grad \Phi_\ep + \Grad \left( \gamma {\tilde \vr}^{\gamma - 1} S_\ep \right) = \ep \mathcal{G}^2_\ep,
\eF
supplemented with the Neumann boundary condition
\bFormula{AA11}
\Grad \Phi_\ep \cdot \vc{n} |_{\partial \Omega} = 0,
\eF
where we have set
\[
\mathcal{G}^1_\ep = \frac{1}{\tilde \vr} \Div \left( \vre
\left( \frac{1 - \Theta_\ep}{\ep} \right) \vue \right),
\]
and
\[
\mathcal{G}^2_\ep = - \mathcal{Q}_{\tilde \vr} \Grad \left[ \frac{(\vre \Theta_\ep)^\gamma  - \gamma {\tilde \vr}^{\gamma - 1}(
\vre \Theta_\ep - \tilde \vr) - {\tilde \vr}^\gamma}{\ep^2} \right] -  \mathcal{Q}_{\tilde \vr} \Div (\vre \vue \otimes \vue ) -
\mathcal{Q}_{\tilde \vr} \Div \tn{S}(\Grad \vue)
\]
\[
+ \ep \mathcal{Q}_{\tilde \vr} \vre \frac{1 - \The}{\ep^2} \Grad F.
\]

\subsection{Analysis of the acoustic system}

Consider the homogeneous problem associated to (\ref{AA9}), (\ref{AA10}), namely
\bFormula{AA12}
\partial_t S + \frac{1}{\tilde \vr} \Div ( \tilde \vr \Grad \Phi) = 0, \
\partial_t \Phi + c(\tilde \vr) S  = 0,\ \Grad \Phi \cdot \vc{n}|_{\partial \Omega} = 0,
\eF
where we have denote
\[
c(\tilde \vr) = {p'(\tilde\vr)}=\gamma {\tilde \vr}^{\gamma - 1},
\]
that can be viewed as a wave equation
\[
\partial^2_{t,t} \Phi - \frac{c (\tilde \vr) }{\tilde \vr} \Div ( \tilde \vr \Grad \Phi) = 0, \ \Grad \Phi \cdot \vc{n}|_{\partial \Omega} = 0
\]
for the potential $\Phi$.

The acoustic propagator
\[
\mathcal{A} (w) = - c(\tilde \vr) \Delta_y w - \frac{c(\tilde \vr)}{\tilde \vr} \partial_z \left( \tilde \vr \partial_z w \right),
\quad {x=(y_1,y_2,z)}, \quad \Grad w \cdot \vc{n}|_{\partial \Omega} = 0,
\]
has been studied by several authors. In particular, DeBi\' evre and Pravica \cite{DebPr1}, \cite{DebPr2} showed that $\mathcal{A}$ can be viewed as a self-adjoint operator
on a weighted $L^2-$space endowed with the scalar product
\[
\left< u, v \right>_H = \intO{ u v \frac{\tilde \vr}{c(\tilde \vr)}}.
\]

Following the strategy of \cite{Fei2011} we derive the desirable dispersive estimates for the solutions of the acoustic system (\ref{AA9}), (\ref{AA10})
from the following version of the celebrated \emph{RAGE theorem} (see Cycon et al. \cite[Theorem 5.8]{CyFrKiSi}):

\bTheorem{RAGE}
Let $H$ be a Hilbert space, $A: {\cal D}(A) \subset H \to H$ a
self-adjoint operator, $C: H \to H$ a compact operator, and $P_c$
the orthogonal projection onto the space of continuity $H_c$ of $A$,
specifically,
\[
H = H_c \oplus {\rm cl}_H \Big\{ {\rm span} \{ w \in H \ | \ w \ \mbox{an eigenvector of} \ A \} \Big\}.
\]

Then
\bFormula{rage1}
\left\| \frac{1}{\tau} \int_0^\tau \exp(-{\rm i} tA ) C P_c \exp(
{\rm i} tA ) \ \dt \right\|_{{\cal L}(H)} \to 0 \ \mbox{as}\ \tau
\to \infty.
\eF
\eT

\subsubsection{Absence of the point spectrum for the wave propagator}
\label{W}

In order to apply Theorem \ref{TRAGE} we have to make sure that the point spectrum of $\mathcal{A}$ is empty. To see this
we use the positive commutator method introduced by DeBi\' evre and Pravica
\cite{DebPr2} in a similar context. More precisely, we compute the commutator $\left[ \mathcal{A}, y \cdot \nabla_y \right]$:
\bFormula{W1}
 \left[ \mathcal{A}, y \cdot \nabla_y \right](w) \equiv \mathcal{A}( y \cdot \nabla_y w) -
y \cdot \nabla_y (\mathcal{A}(w))
\eF
\[
= - c(\tilde \vr) \left[ \Delta_y ( y \cdot \nabla_y w) - y \cdot \nabla_y (\Delta_y w) \right] = -2 c(\tilde \vr)  \Delta_y w.
\]
If
\[
\mathcal{A}(w) = \lambda w,
\]
then, {since $\mathcal{A}$ is self-adjoint}, we get
\[
\int_\Omega {w} \left[ \mathcal{A} , y \cdot \nabla_y \right]{{w}} \ {\rm d}y \ {\rm d}z = \int_\Omega \left( \mathcal{A}(w) y \cdot \nabla_y w - \lambda {w} y \cdot \nabla_y w \right) \ {\rm d}y \ {\rm d}z = 0,
\]
while, in accordance with (\ref{W1}),
\[
\int_\Omega {w} \left[ \mathcal{A} , y \cdot \nabla_y \right]{{w}} \ {\rm d}y \ {\rm d}z = 2 \int_\Omega c(\tilde \vr) |\nabla_y w |^2.
\]
Thus any possible eigenfunction $w$ must be constant with respect to the horizontal component; whence $w \equiv 0$.

\subsubsection{Decay of acoustic waves}

Following step by step \cite[Section 5]{Fei2011} we rewrite the acoustic system (\ref{AA9}), (\ref{AA10}) in the abstract form
\bFormula{AA20}
\ep \partial_t \left< Z_\ep , \phi \right>_H  + \left< \sqrt{ \mathcal{A}}[\Phi_\ep], \sqrt{\mathcal{A}}[\phi] \right>_H = \ep \left< \mathcal{F}^1_\ep , \phi \right>_H, \eF
for any smooth $\phi \in \DC(\Ov{\Omega}) ,\  \Grad \phi \cdot \vc{n}|_{\partial \Omega} = 0$,
\bFormula{AA21}
\ep \partial_t \left< \Phi_\ep, \phi \right>_H + \left< Z_\ep, \phi  \right>_H = \ep \left< \mathcal{F}^2_\ep, \phi \right>_H,
\ \mbox{for any smooth} \ \phi \in \DC(\Ov{\Omega}),
\eF
where we have set $Z_\ep = c(\tilde \vr) S_\ep$. The quantities $\mathcal{F}^1_\ep$, $\mathcal{F}^2_\ep$ are computed in terms of
$\mathcal{G}^1_\ep$, $\mathcal{G}^2_\ep$. After a bit tedious but straightforward computation for which we refer to \cite[Section 5]{Fei2011}, it can be shown that
\[
\mathcal{F}^1_\ep = H^1 (\mathcal{A}) [h^1_\ep], \ \mathcal{F}^2_\ep = H^2 (\mathcal{A} ) [h^2_\ep],
\]
where
\[
\{ h^1_\ep \}_{\ep > 0}, \ \{ h^2_\ep \}_{\ep > 0} \ \mbox{are bounded in} \ L^2((0,T) \times \Omega),
\]
and $H^1, H^2 : (0, \infty) \to \R$ are suitable smooth functions that may become singular at $0$ and $\infty$.

Now, exactly as in \cite[Section 5.2]{Fei2011}, we apply Theorem \ref{TRAGE} in the situation
\[
A = \sqrt{\mathcal{A}},\ C = \chi^2 G(\mathcal{A}), \ \chi \in \DC(\Ov{\Omega}), \chi \geq 0, \ G \in \DC(0,\infty), \ G \geq 0,
\]
to deduce that
\bFormula{AA22}
\left\| G(\mathcal{A}) [ \Phi_\ep ] \right\|_{L^2((0,T) \times K) } \to 0 \ \mbox{as}\ \ep \to 0,
\eF
\bFormula{AA23}
\left\| G(\mathcal{A}) [ Z_\ep ] \right\|_{L^2((0,T) \times K) } \to 0 \ \mbox{as}\ \ep \to 0
\eF
for any compact $K \subset \Ov{\Omega}$ and any $G \in \DC(0, \infty)$, cf. \cite[Formula (5.20)]{Fei2011}. Moreover, seeing that $\Phi_\ep = \mathcal{Q}_{\tilde \vr}(\vre \vue)$ we may use estimate (\ref{wc1}) to deduce from (\ref{AA22}) that
\bFormula{AA24}
\mathcal{Q}_{\tilde \vr}(\vre \vue) \to 0 \ \mbox{in}\ L^2_{\rm loc}((0,T) \times K; \R^3).
\eF

Relations (\ref{Aa1}), (\ref{AA24}) imply the convergence claimed in (\ref{mw2}), and, accordingly, (\ref{final2}). In order to complete
the proof of Theorem \ref{Tw2}, it remains to show (\ref{final1}). To this end, we will use the convergence result (\ref{AA23}) together with some uniform bounds derived in the following section.

\subsection{Uniform pressure estimates}
\label{UPE}

Our ultimate goal is to show (\ref{final1}). Obviously, using the uniform bounds (\ref{B8}--\ref{B15}), the convergence (\ref{AA23}) implies that
\bFormula{AAA1}
\frac{ \vre \The - \tilde \vr}{\ep \tilde \vr} \to 0 \ \mbox{in}\ L^2(0,T; L^2_{\rm weak} (K)) \ \mbox{for any compact}\
K \subset \Ov{\Omega}.
\eF
Consequently, in order to deduce (\ref{final1}) we have to establish:
\[
\mbox{{\bf (i)}
spatial compactness of}\ {\frac{ \vre \The - \tilde \vr}{\ep \tilde \vr}},
\]
\[
\mbox{
{\bf (ii)}
better integrability of the expression}\
\frac{ (\vre \Theta_\ep)^\gamma  - \gamma {\tilde \vr}^{\gamma - 1}(
\vre \Theta_\ep - \tilde \vr) - {\tilde \vr}^\gamma }{\ep^2}.
\]

As for the spatial compactness claimed in {\bf (i)}, we refer to Masmoudi \cite[Proposition 4.1]{MAS2}. Specifically, denoting
$\{ \omega_\delta \}_{\delta > 0}$ the family of regularizing kernels, it can be shown that
\bFormula{smooth}
\left\| \omega_\delta * \left(\frac{ \vre \The - \tilde \vr}{\ep }\right) - \left(\frac{ \vre \The - \tilde \vr}{\ep}\right)
\right\|_{L^p(0,T; L^2(K))} \to 0 \ \mbox{as}\ \delta \to 0,\ 1 \leq p < \infty,
\eF
uniformly for $\ep > 0$ for any compact $K \subset \Omega$, which, together with (\ref{AAA1}) gives rise to
\bFormula{AAA2}
\frac{ \vre \The - \tilde \vr}{\ep} \to 0 \ \mbox{(strongly)  in}\ L^2((0,T) \times K)) \ \mbox{for any compact}\
K \subset \Ov{\Omega}.
\eF

As for {\bf (ii)}, we can proceed as in \cite[Section 6, formula (6.8)]{Fei2011} to show that
\bFormula{AAA3}
\int_0^T \int_{ \{ \vre \The < \tilde \vr/2 \} \cap K } (\vre \The)^{\gamma + 1} \ \dxdt + \int_0^T \int_{ \{ \vre \The > 2 \tilde \vr \} \cap K } (\vre \The)^{\gamma + 1} \ \dxdt \leq \ep^2 c(K)
\eF
for any compact $K \subset \Ov{\Omega}$.


We can now combine (\ref{AAA2}), (\ref{AAA3}) to show  (\ref{final1}) on compact subsets of $\Omega$. This can be shown by repetition of the argument from \cite[Section 6,  formulas (6.5) and (6.9)]{Fei2011} when $\vr_\ep$ is replaced by $\vr_\ep\Theta_\ep$.
 We have proved Theorem \ref{Tw2}.

\section{Concluding remarks}
\label{C}

The rather inconvenient restriction {$\gamma> 3$}  in Theorem \ref{Tw2}
is purely technical and has been effectively used only in Section \ref{UPE}, more specifically, only in the proof of spatial compactness
(\ref{smooth}). It could have been relaxed should we have more information on integrability of
$\vre$. An alternative approach to the problem would be to use the relative energy inequality, similarly to \cite{FeJiNo1}. Such an approach, however, would require Strichartz type space-times estimates for the operator $\mathcal{A}$ in $\Omega$. To best of our knowledge, validity of these estimates is not known and maybe even not true in general, at least for the unbounded slab $\Omega = \R^2 \times (0,1)$.

Our final remark concerns \emph{existence} of the weak solutions enjoying the regularity required in the present paper, in particular satisfying
(\ref{w5}) or, alternatively, (\ref{w5a}). As already mentioned in the introduction, an existence result in the variables $[\vr, \vu, Z =
\vr \Theta]$ has been established in \cite{MMNP}. Next, we claim the following result that may be seen as a straightforward consequence
of the theory developed by DiPerna and Lions \cite{DL}.

\bLemma{CC}
Suppose
\bFormula{CC1}
\vr \geq 0, \ \vr \in L^2((0,T) \times \Omega), \ \underline{\Theta} \vr \leq Z \leq \Ov{\Theta} \vr,\ 0 < \underline{\Theta} \leq 1 \leq \Ov{\Theta},
\eF
and
\bFormula{CC2}
\partial_t \vr + \Div (\vr \vu) = 0, \ \partial_t Z + \Div (Z \vu) = 0 \ \mbox{in}\ (0,T) \times \R^3,
\eF
where the velocity field belongs to the class
\bFormula{CC3}
\vu \in L^2(0,T; W^{1,2}(\R^3;\R^3)).
\eF

Then
\[
\Theta = \left\{ \begin{array}{l} Z / \vr  \ \mbox{for}\ \vr > 0, \\ \\
1 \ \mbox{for}\ \vr = 0 \end{array} \right. \in L^\infty((0,T) \times \Omega)
\]
satisfies the transport equation
\bFormula{CC4}
\partial_t \Theta + \vu \cdot \Grad \Theta = 0.
\eF

\eL

\bProof

Using the regularizing procedure of DiPerna and Lions \cite{DL} we obtain
\bFormula{CC5}
\partial_t \vr_\delta + \vu \cdot \Grad \vr_\delta + \vr_\delta \Div \vu = r^1_\delta,
\eF
\bFormula{CC6}
\partial_t Z_\delta + \vu \cdot \Grad Z_\delta + Z_\delta \Div \vu = r^2_{\delta},
\eF
where, in view of (\ref{CC1}), (\ref{CC3}) and Friedrich's lemma,
\[
r^1_\delta, \ r^2_\delta \to 0 \ \mbox{in}\ L^1((0,T) \times \R^3) \ \mbox{as}\ \delta \to 0.
\]

Next, we multiply (\ref{CC5}) by $- (Z_\delta + \lambda) /(\vr_\delta + \lambda)^2$, (\ref{CC6}) by $(\vr_\delta + \lambda)^{-1}$, respectively, where
$\lambda > 0$ is a positive constant. Summing up the resulting expressions we deduce, after a straightforward manipulation,
\[
\partial_t \left( \frac{ Z_\delta + \lambda }{\vr_\delta + \lambda } \right) + \Div \left[ \left( \frac{ Z_\delta + \lambda }{\vr_\delta + \lambda } \right) \vu \right]
- \left[ \frac{( Z_\delta + \lambda) \vr_\delta }{(\vr_\delta + \lambda)^2} + \frac{\lambda}{ \vr_\delta + \lambda }\right] \Div \vu
\]
\[
= - r^1_\delta \frac{ Z_\delta + \lambda }{(\vr_\delta + \lambda)^2 } + r^2_\delta \frac{1}{\vr_\delta + \lambda}.
\]
By virtue of (\ref{CC1}), we may let first $\delta \to 0$ and then $\lambda \to 0$ to recover the desired equation (\ref{CC4}).

\qed

The solutions constructed in \cite{MMNP} enjoy the properties (\ref{CC1}), (\ref{CC2}), in particular, Lemma \ref{LCC} yields (\ref{w5a})
at least for the values of the adiabatic exponent $\gamma \geq 2$. In order to prove existence for the full range $\gamma > 3/2$, one has
to introduce the ``artificial pressure'' proportional to $\delta Z^2$ similarly to \cite[Chapter 3]{FeNo6} and then perform the limit
$\delta \to 0$ adapting the technique of Mich\' alek \cite{Mich}.

\def\cprime{$'$} \def\ocirc#1{\ifmmode\setbox0=\hbox{$#1$}\dimen0=\ht0
  \advance\dimen0 by1pt\rlap{\hbox to\wd0{\hss\raise\dimen0
  \hbox{\hskip.2em$\scriptscriptstyle\circ$}\hss}}#1\else {\accent"17 #1}\fi}

\end{document}